# Algebraic Constructions of Efficient Broadcast Networks

Michael J. Dinneen and Michael R. Fellows
Department of Computer Science
University of Victoria
Victoria, B.C. Canada V8W 3P6

Vance Faber
Los Alamos National Laboratory
Los Alamos, New Mexico 87545, U.S.A.

*Abstract.* Cayley graph techniques are introduced for the problem of constructing networks having the maximum possible number of nodes, among networks that satisfy prescribed bounds on the parameters maximum node degree and broadcast diameter. The broadcast diameter of a network is the maximum time required for a message originating at a node of the network to be relayed to all other nodes, under the restriction that in a single time step any node can communicate with only one neighboring node. For many parameter values these algebraic methods yield the largest known constructions, improving on previous graph-theoretic approaches. It has previously been shown that hypercubes are optimal for degree $k$ and broadcast diameter $k$. A construction employing dihedral groups is shown to be optimal for degree $k$ and broadcast diameter $k+1$.

## 1. Introduction.

The problem of designing efficient networks arises in many different contexts, including parallel processing, communication networks and security systems. Several design constraints arise in practice, such as bounds on the maximum degree of network vertices, planarity and symmetry properties. A variety of resource efficiencies may be the design objective, such as the numbers of vertices and edges, communication times, and fault tolerance. Good constructions for small parameter values are important, because of engineering applications for network designs of small and moderate size.

An instance of this general research program is the extensively studied problem of finding constructions of graphs having the maximum possible number of vertices for given bounds on maximum vertex degree and diameter. For example, the Petersen graph of order 10 is the unique largest possible graph of maximum degree 3 and diameter 2. See figure 1(a) below. The Petersen graph is vertex symmetric. The labeling in the figure indicates the distance to the vertex labeled 0. For recent surveys concerning this *degree/diameter* problem see [BDQ,Ch].

Previous work of the authors (and others) has shown that algebraic techniques can be applied powerfully to the degree/diameter problem [BDV,CCD,CCDFFLMMS,CCSW]. Cayley graphs now yield most of the largest known constructions for small parameter values, surpassing many earlier results based on graph composition operators (see [BDQ]). In this paper we introduce algebraic techniques for the analogous network design problem concerning broadcasting. Specifically, we address the problem of finding constructions of graphs having the maximum possible number of vertices for given bounds on maximum vertex degree and broadcast time.

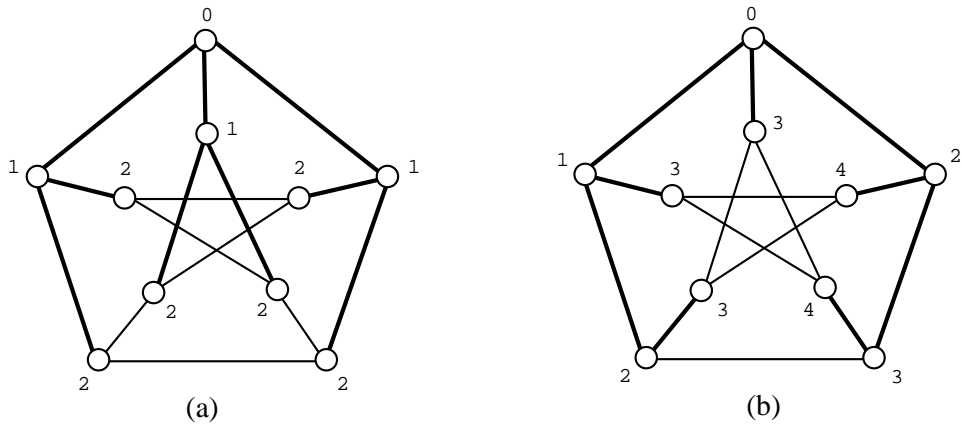

Figure 1. Petersen graph showing (a) diameter 2 (b) broadcast time 4.

Broadcasting concerns the dissemination of a message originating at one node of a network to all other nodes. This task is accomplished by placing a series of calls over the communication lines of the network between neighboring nodes, where each call requires a unit of time and a call can involve only two nodes. A node can participate in only one call per unit of time. Thus at time $t$ it is possible for at most $2^t$ nodes to have received the message, in any network.

Figure 1(b) shows that the broadcast time of the Petersen graph is 4. The edges belonging to the tree of calls of a time 4 broadcast are indicated; the label of a node shows the time at which the message is received. Figure 1 invites the following simile. The *diameter* is the time required to disseminate a message by shouting, while the *broadcast time* measures the time required if dissemination is by telephoning.

For bibliographic surveys of the extensive literature on broadcasting in networks and closely related problems see [HHL,LS].

In section 2 we introduce terminology and notation, and describe our algebraic approach to broadcast network construction. In section 3 we prove, using dihedral groups, a new infinite optimal family.

## 2. Preliminaries and An Algebraic Approach.

Graphs in this paper are simple and undirected. We use the terms *network* and *graph* interchangeably, and similarly *node* and *vertex*. Let $G = (V, E)$ be a connected graph and let $u$ be a vertex of $G$. The *broadcast time of vertex* $u$, $b(u)$, is the minimum number of time units required to complete broadcasting of a message originating at vertex $u$. The *broadcast time of* $G$ is the maximum broadcast time of any vertex $u$ in $G$, $b(G) = \max\{b(u) \mid u \in V\}$.

We study the function $B(\Delta, t)$ defined to be the maximum possible order of a graph $G$ of maximum degree $\Delta$ and broadcast time at most $t$. A graph of maximum degree $\Delta$ and broadcast time $t$ is termed an *optimal* $(\Delta, t)$-broadcast network if it has order $B(\Delta, t)$.

It is straightforward to derive a recurrence relation for a function $M(\Delta, t)$ that is an analogue of the Moore bound for the degree/diameter problem. We record this well-known bound in the following proposition.

*Proposition 1.*  $B(\Delta, t) \leq M(\Delta, t)$ where
(1) $f(\Delta, 0) = 1$ for all $\Delta$,
(2) $f(\Delta, t) = 1 + \Sigma_{i=1}^{\min(\Delta, t)} f(\Delta, t-i)$ and
(3) $M(\Delta, t) = 2 \cdot f(\Delta - 1, t - 1)$  □

It is easy to observe that if $G$ is a graph with broadcast time $t$ then the Cartesian product of $G$ and $K_2$ has broadcast time at most $t + 1$. This has the following useful and well-known consequence.

*Proposition 2.*  $B(\Delta + 1, t + 1) \geq 2 \cdot B(\Delta, t)$  □

The $r$-dimensional hypercube $Q_r$ is isomorphic to the Cartesian product of $r$ copies of $K_2$. Since any graph with broadcast time $t$ has at most $2^t$ vertices, repeated application of Prop. 2 yields the following well-known fact.

*Proposition 3.*  $Q_\Delta$ is an optimal $(\Delta, \Delta)$-broadcast network.  □

The hypercubes provide a simple example of a general method of network construction based on groups. Let $A$ be a group and let $S \subseteq A$ be a set of generators of $A$ that is closed under taking of inverses (if $a \in A$ then also $a^{-1} \in A$). The *Cayley graph* $(A, S)$ is the graph with vertex set the elements of $A$, having an edge between the pair of elements $a, b$ of $A$ if and only if $as = b$ for some $s \in S$. The hypercubes are thus Cayley graphs on the $Z_2$ vector spaces, with generator sets consisting of the standard basis elements.

Cayley graphs are *vertex symmetric*. That is, there is an automorphism of the graph taking any vertex to any other. Symmetry may be an important property of networks in applications such as parallel processing [ABR]. One of the advantages of a Cayley graph is that it suffices to find a broadcast scheme requiring time $t$ for the identity node of the network. This can then be translated by group multiplication to provide a broadcast schedule for any other node originating a message.

Suppose $(A, S)$ is a Cayley graph with generators $S = \{s_1, ..., s_k\}$. One broadcasting method is to consider the indexing of the generators as the order in which calls should be made to neighboring vertices in a broadcast. Given a group with an easily computable multiplication, it is relatively simple to compute the broadcast time of the Cayley graph using this scheme by starting with the identity element and proceeding until all elements of the group have been generated. Experimental computing using this technique led us to discover the new infinite family of optimal constructions described in the next section.

Another method is to choose $t$ different permutations $\pi_i$ of the generator set $S$. A node receiving the message at time $i$ places calls to its neighbors in the order indicated (by multiplication) by the sequence of generators given by $\pi_i$. Experimental computing using this method (and further generalizations) has provided many new constructions that are the largest known for some parameter values. These are indicated by an asterisk in Table 1 below. Details concerning these constructions can be found in [Di].

*Example.* The group $Z_{12}$ acts on the group $Z_{13}$ by multiplication. Thus we may form the semi-direct product $Z_{12} \times_\alpha Z_{13}$ where $\alpha$ denotes this action, $\alpha_k(x) = 2^k \cdot x \bmod 13$. The generators: (7,1), (5,7) and (6,0) give us a Cayley graph of order 156, degree 3 and broadcast time 10, which is presently the largest known (3,10)-broadcast network.

Table 1 shows the orders of the largest known broadcast graphs for small parameter values. Optimal values are displayed in bold. Entries below the diagonal are omitted as they trivially follow from the diagonal entries. Table 2 shows the best known upper bounds for these parameter values.

| $\Delta \setminus t$ | 2 | 3 | 4 | 5 | 6 | 7 | 8 | 9 | 10 |
|---|---|---|---|---|---|---|---|---|---|
| 2 | **4** | **6** | **8** | **10** | **12** | **14** | **16** | **18** | **20** |
| 3 |   | **8** | **14** | **24** | **40** | 60* | 84* | 126* | 156* |
| 4 |   |   | **16** | **30** | **56** | 90* | 148* | 253* | 272* |
| 5 |   |   |   | **32** | **62** | 108* | 186* | 336* | 506 |
| 6 |   |   |   |   | **64** | **126** | 220* | 390* | 750* |
| 7 |   |   |   |   |   | **128** | **254** | 440 | 816* |
| 8 |   |   |   |   |   |   | **256** | **510** | 880* |
| 9 |   |   |   |   |   |   |   | **512** | **1022** |

Table 1. Largest known broadcast networks.

| $\Delta \setminus t$ | 2 | 3 | 4 | 5 | 6 | 7 | 8 | 9 | 10 |
|---|---|---|---|---|---|---|---|---|---|
| 2 | 4 | 6 | 8 | 10 | 12 | 14 | 16 | 18 | 20 |
| 3 | 4 | 8 | 14 | 24 | 40 | 66 | 108 | 176 | 286 |
| 4 | 4 | 8 | 16 | 30 | 56 | 104 | 192 | 354 | 652 |
| 5 | 4 | 8 | 16 | 32 | 62 | 120 | 232 | 448 | 864 |
| 6 | 4 | 8 | 16 | 32 | 64 | 126 | 248 | 488 | 960 |
| 7 | 4 | 8 | 16 | 32 | 64 | 128 | 254 | 504 | 1000 |
| 8 | 4 | 8 | 16 | 32 | 64 | 128 | 256 | 510 | 1016 |
| 9 | 4 | 8 | 16 | 32 | 64 | 128 | 256 | 512 | 1022 |
| 10 | 4 | 8 | 16 | 32 | 64 | 128 | 256 | 512 | 1024 |

Table 2. Upper bounds on $B(\Delta, t)$.

## 3. A New Infinite Family of Optimal Networks.

Previously, the hypercubes were the only infinite family of graphs known to be optimal. Experimental computing with the Cayley graph techniques described in the last section led to the formulation of our main theorem. The first few graphs in our infinite family were discovered independently in [BHLP1, BHLP2] by graph theoretic methods.

*Theorem 1.* The Cayley graphs from the dihedral groups $D_{n_\Delta}(n_\Delta = 2^\Delta - 1)$ with generators $w, wx^1, wx^3, \ldots, wx^{2^{\Delta-1}-1}$ where $w^2 = e$, $wxw^{-1} = x^{-1}$ and $x^{2^\Delta-1} = e$, are a family of optimal $(\Delta, \Delta + 1)$-broadcast networks.

*Proof.* First note that each of the generators is an involution. We describe a broadcast scheme for a message originating at the identity element. For the first $\Delta$ time steps we follow the rule that each element $v$ of the group that knows the message, communicates to the neighboring element $v \cdot wx^{2^{k-1}-1}$ in the Cayley graph at time step $k$. Thus, letting $T_i = \{v \mid \text{vertex } v \text{ knows the message at time } i \}$, we have

$$
\begin{aligned}
T_0 &= \{e = w^0 = x^0\} \\
T_1 &= T_0 \cup \{w\} \\
T_2 &= T_1 \cup \{wx^1, x^1\} \\
T_3 &= T_2 \cup \{wx^3, x^3, wx^{-1+3}, x^{-1+3}\} \\
&= \bigcup_{i=0}^{3} \{wx^i, x^i\} \\
&\vdots \\
T_k &= T_{k-1} \cup T_{k-1} \cdot wx^{2^{k-1}-1} \\
&= \bigcup_{i=0}^{2^{k-2}-1} \{wx^i, x^i\} \cup \{wx^{2^{k-1}-1-i}, x^{2^{k-1}-1-i}\} \\
&= \bigcup_{i=0}^{2^{k-1}-1} \{wx^i, x^i\} \\
&\vdots \\
T_\Delta &= T_{\Delta-1} \cup T_{\Delta-1} \cdot wx^{2^{\Delta-1}-1} \\
&= \bigcup_{i=0}^{2^{\Delta-1}-1} \{wx^i, x^i\}
\end{aligned}
$$

For the last step, every vertex $v$ in $T_\Delta$ transmits the message to the neighbor given by multiplication with $w$, so that

$$
\begin{aligned}
T_{\Delta+1} &= T_\Delta \cup T_\Delta \cdot w \\
&= \bigcup_{i=0}^{2^{\Delta-1}-1} \{wx^i, x^i\} \cup \bigcup_{i=1}^{2^{\Delta-1}-1} \{wx^i w, x^i w\} \\
&= \bigcup_{i=0}^{2^{\Delta-1}-1} \{wx^i, x^i\} \cup \bigcup_{i=1}^{2^{\Delta-1}-1} \{wx^{n_\Delta-i}, x^{n_\Delta-i}\} \\
&= \bigcup_{i=0}^{2^\Delta-1} \{wx^i, x^i\} = D_{n_\Delta}.
\end{aligned}
$$

We see that after $\Delta + 1$ time steps all of the vertices have received the message. The networks are optimal since

$$|D_{n_\Delta}| = 2(2^\Delta - 1) = 2^{\Delta+1} - 2$$

is also the upper bound on $B(\Delta, \Delta + 1)$. $\square$

## 4. Conclusions and Open Problems.

Our main contribution in this brief paper is the demonstration that elementary Cayley graph techniques can be applied to the problem of designing efficient broadcast networks, giving improvements over previous graph theoretic methods. Further computational exploration will undoubtedly reveal further record-breaking constructions for small parameter values, and perhaps suggest additional infinite families of optimal constructions. It would be interesting to see if similar techniques can be fruitfully applied to related network design problems that are as yet unexplored from an algebraic perspective, such as broadcasting in directed networks [LP] and gossiping [HHL].